\documentclass[twoside]{article}
\usepackage[papersize={5.5in,8.27654in},total={4.25in,7.25in},%
  includehead,nofoot,centering]{geometry}
\usepackage{amsmath,txfonts}
\usepackage[T1]{fontenc}
\usepackage[colorlinks=true,allcolors=blue,bookmarks=false,
  pdfauthor={Gavin R. Putland},
  pdftitle={Generalized Gregorian quadrature, including end-corrected weights for the midpoint rule},
  pdfkeywords={Gregory’s method, quadrature, midpoint rule, trapezoidal rule, Newton-Cotes, Simpson’s rule},
  pdfsubject={Numerical Analysis, Numerical Integration}]{hyperref}
\usepackage[final,kerning,spacing]{microtype}
  \SetExtraKerning{encoding={*}}{1={-35,-45}, / ={50,50},
    \textemdash={50,50}, \textendash={30,30}, \textquotedblleft={-30,}}
  \SetExtraSpacing{encoding={*}}{A={-80,,}, f={80,,}, r={-80,,}, v={-80,,},
    w ={-70,,}, y={-60,,}, ' = {-200,,}, \textquotedblright={-80,,}}

\title{Generalized Gregorian quadrature, including end-corrected weights for the midpoint rule}
\author{Gavin R.{\small~}Putland~\!\thanks{\,\small Royal Melbourne Institute of Technology, Australia.~ Gmail name: grputland.\\Copyright license: \href{https://creativecommons.org/licenses/by/4.0/legalcode}{Creative Commons Attribution 4.0 International}.}}
\date{\vspace{-1ex}\normalsize 17 December 2025} 

\makeatletter
\def\ps@plain{
  \let\@oddfoot=\@empty
  \def\@oddhead{\hfill\normalsize\sf\thepage}
  }
\def\ps@headings{
  \def\sectionmark##1{\markright{\S\thesection.~ ##1}}
  \let\@evenfoot=\@empty  \let\@oddfoot=\@empty
  \def\@evenhead{\footnotesize\sf\underline{\makebox[\textwidth]{\normalsize\sf\thepage\hfill\footnotesize\sf Putland,\itshape\, Generalized Gregorian quadrature\,\ldots}}}
  \def\@oddhead{\footnotesize\sf\underline{\makebox[\textwidth]{\rightmark\hfill\normalsize\sf\thepage}}}
  }
\makeatother

\begin{document}

\sloppy

\maketitle

\subsection*{Abstract}

A class of numerical quadrature rules is derived, with equally-spaced nodes, and unit weights except at a few points at each end of the series, for which "corrections" (not using any further information about the integrand) are added to the unit weights. If the correction sequences overlap, the effects are additive. A fundamental parameter ("alpha") in the derivation is the distance from the endpoint of the range of integration to the first node, measured inward in step-lengths. Setting alpha to 1/2 yields a set of corrected composite midpoint rules. Setting alpha=0 yields Gregory's closed Newton-Cotes-like rules, including (for sufficient overlap) the standard closed Newton-Cotes rules (trapezoidal rule, "1/3 Simpson rule", "3/8 Simpson rule", "Boole's rule", etc.). Setting alpha=1 yields open N-C-like rules, again including the standard ones. A negative alpha means that the integrand is sampled outside the range of integration; suitably chosen negative values yield centered finite-difference end-corrections for the trapezoidal rule and the midpoint rule. One can even have different values of alpha at the two ends, yielding, inter alia, Adams-Bashforth and Adams-Moulton weights. Thus the title could have been "Unified derivation of equispaced quadrature rules".

\clearpage

\pagestyle{headings}

\section{Framing the problem}

We seek a numerical quadrature rule of the form
\begin{equation}\label{e-cd}
  \int_{-\alpha h}^{(n\text{+}\beta)h} \!\!f(t)\,dt \,\approx~
    h\sum_{i\text{=0}}^n f(ih) \,+\,
    h\sum_{i\text{=0}}^m c_i ~\!f(ih) \,+\,
    h\sum_{i\text{=0}}^m d_i ~\!f\big((n\!-\!i)h\big) \,,
\end{equation}
where the coefficients $c_i$ and $d_i$ are independent of $n$ and~$h$
(the step size), but may depend on $\alpha$ and~$\beta$.\, Such a rule
would have the following advantages:
\begin{itemize}
\item While sampling the integrand at equally-spaced abscissae
  (\textbf{nodes}) would be convenient\textemdash and might even be
  dictated by the data\textemdash the parameters $\alpha$ and~$\beta$
  would allow the \textbf{terminals} (``limits'' of integration) to be
  at \emph{or~between} nodes.
\item Particular values of the parameters would yield useful special
  cases. For~${\alpha\!=\!\beta\!=\!0}$, we~would get ``closed
  Newton-Cotes--like rules''\!, for which the outermost nodes coincide
  with the terminals. For~${\alpha\!=\!\beta\!=\!1}$, we~would get
  ``\emph{open} N-C--like rules''\!, for which the outermost nodes are
  one step in from the terminals. (N-C--like rules are familiar; but
  reproducing the most famous examples would serve as a sanity check
  on our method.) Negative parameters (or $m$ greater than~$n$) would
  yield rules with nodes outside the range of integration.
  For~${\alpha\!=\!\beta\!=\!\text{1/2}}$, the range of integration
  would be divided into $n\text{+1}$ steps of length~$h$, with the
  nodes centered in the steps, yielding a modified \textit{midpoint
    rule}\textemdash as advertised in the title.
\item For\, $n>2m\!+\!2$,\, the right side of\,\eqref{e-cd} divided
  by~$h$\, would become
  \begin{equation}
    \sum_{i=m\text{+1}}^{n-m-1} f(ih)
    \,+ \sum_{i\text{=0}}^m (\text{1+}c_i)~\!f(ih)
    \,+ \sum_{i\text{=0}}^m (\text{1+}d_i)~\!f\big((n\!-\!i)h\big) \,,
  \end{equation}
  which is a rule with \emph{unit weights} except at $m\text{+1}$
  nodes at each end, where the coefficients $c_i$ and $d_i$ can be
  described as \emph{corrections} to unit weights, or
  \emph{differences} from unit weights (hence the symbols). The
  sequence of unit weights in the interior would have no cycles with
  periods longer than the step size, minimizing risk of bias due to
  cycles in the weights{\tiny\,} interacting with oscillations in the
  integrand,\footnote{\,\small In signal-processing terms, a cycle in
  the weights is tantamount to a reduction in the sampling rate.}
  eliminating the need for the number of steps to be a multiple of any
  cycle length, and expediting the task of entering the weights into a
  spreadsheet (``Fill down''!). For\, $n\le 2m\!+\!2$,\, the
  correction sequences would overlap and, according to~\eqref{e-cd},
  the corrections would be additive.
\item Hence, for~${\alpha\!=\!\beta\!=\!\text{1/2}}$, we would get end
  corrections for the composite midpoint rule (rectangle rule). And
  for~${\alpha\!=\!\beta\!=\!0}$, we would get end corrections for the
  composite trapezoidal rule, taking unit weights as the base
  case. These corrections would not need further information about the
  integrand, but would merely adjust the weights\textemdash unlike the
  standard ``corrected trapezoidal rule''\!, which requires
  derivatives at the terminals.
\item By using two successive values of~$m$, we could compare two
  estimates of the integral from the \emph{same ordinates} for the
  purpose of error control
  (\textit{cf}.~Runge--Kutta--Fehlberg\,/\,Runge--Kutta--\!Verner
  methods for ordinary differential equations). This procedure, unlike
  comparing estimates from the same ordinates using N--C methods with
  different orders or different cycle lengths, would still not
  restrict the number of steps or introduce cycles in the interior
  weights.
\end{itemize}

If $f(t)$ is a polynomial of degree~$m$ in~$t$,\, each side
of~\eqref{e-cd} is a polynomial of degree $m\text{+1}$ in~$h$, with no
constant term. For a more general~$f(t)$, the expansion of each side
of~\eqref{e-cd} in powers of~$h$ will still have no constant term. So,
making rule~\eqref{e-cd} exact for polynomials of degree~$m$ is a
matter of matching the coefficients of~$h^{k\text{+1}}$ for
${k\!\in\![0\,\text{..}\,m]}$ with a general\,$f$,\, \emph{for
all}~$n$; if this is done, the error will generally
be~$O(h^{m\text{+2}})$.

\section{Existence of the solution}

But, for given $\alpha$ and $\beta$, why should there exist constants
$c_i$ and $d_i$ such that \eqref{e-cd} is exact for all $n$ and~$h$,
for polynomials $f(t)$ of degree up to~$m$? To~answer this, let~us
first split the interval of integration:
\begin{equation}\label{e-split}
  \int_{-\alpha h}^{(n\text{+}\beta)h} \!\!f(t)\,dt ~=
    \int_{-\alpha h}^0 \!f(t)\,dt
    \,+ \int_0^{nh} \!\!f(t)\,dt
    \,+ \int_{nh}^{(n\text{+}\beta)h} \!\!f(t)\,dt \,.
\end{equation}
For the first term on the right, we~replace $f(t)$ by its Taylor
series about ${t\!=\!0}$ (which terminates at the $m^{\text{th}}$
power), and integrate term-by-term, obtaining
\begin{equation}\label{e-alpha}
  \int_{-\alpha h}^0 \!f(t)\,dt \,=\,
    h\sum_{k\text{=0}}^m
      \frac{(-\alpha)^{k\text{+1}}}{(k\text{+1})!} ~\!h^k f^{(k)}(0) \,.
\end{equation}
For the last term we do likewise except that the Taylor series is
about ${t\!=\!nh}$:
\begin{equation}\label{e-beta}
  \int_{nh}^{(n\text{+}\beta)h} \!\!f(t)\,dt \,=\,
    h\sum_{k\text{=0}}^m
    \frac{\beta^{\text{\,}k\text{+1}}}{(k\text{+1})!} ~\!h^k f^{(k)}(nh) \,.
\end{equation}
For the middle term, we know from the Euler--Maclaurin
series~\cite[p.~\!167]{fornberg-21} that
\begin{equation}\label{e-eulermac}\begin{split}
  \int_0^{nh} \!\!f(t)\,dt \,=~\,&
    h\sum_{i\text{=0}}^n f(ih)
    - h\Big[\tfrac{\,1\,}{\,2\,}f(0) + \tfrac{\,1\,}{\,2\,}f(nh)\Big] \\
  & +\text{\,} h\sum_{k\text{=1}}^m a_k ~\!h^k \Big[f^{(k)}(0) - f^{(k)}(nh)\Big] \,,
\end{split}\end{equation}
where the last sum terminates at the degree~$m$ of the polynomial (and
is taken as an empty sum if ${m\!=\!0}$), and the coefficients $a_k$
are constants whose details need not concern us here (except to
acknowledge in passing that ${a_k\!=\!0}$ for positive
\emph{even}~$k$). Now the sum of the right-hand sides of
eqs.\,\eqref{e-alpha} to~\eqref{e-eulermac} is of the form of the
right-hand side of~\eqref{e-cd}, because:
\begin{itemize}
\item The first sum on the right in~\eqref{e-eulermac} is the same as
  in~\eqref{e-cd};
\item In the next term in~\eqref{e-eulermac}, the factor in square
  brackets is a weighted sum of $f(0)$ and~$f(nh)$;
\item $f^{(k)}(0)$ in \eqref{e-alpha} and~\eqref{e-eulermac} is given
  exactly as a weighted sum of the ordinates $f(kh)$ for
  ${k\!\in\![0\,\text{..}\,m]}$, because $f$ itself, being a
  polynomial of degree~$m$, is given exactly as a weighted sum of the
  same $m\text{+1}$ ordinates; and in~order to be dimensionally
  correct, the former weighted sum must have a common factor $h^{-k}$,
  which cancels with~$h^k$;
\item Similarly, $f^{(k)}(nh)$ in \eqref{e-beta}
  and~\eqref{e-eulermac} is given exactly as a weighted sum of the
  ordinates ${f\big((n\!-\!k)h\big)}$ for
  ${k\!\in\![0\,\text{..}\,m]}$, and has a common factor that cancels
  with~$h^k$; and
\item The weights in the aforesaid weighted sums are subsumed under
  $c_i$ and $d_i$ in~\eqref{e-cd}.
\end{itemize}

That explains the form of~\eqref{e-cd} and the conditions under which
it can be made exact. But there are other implications.
In~eqs.\,\eqref{e-alpha} to~\eqref{e-eulermac},\, $\alpha$ appears
only in~\eqref{e-alpha}, where it is related not to $f^{(k)}(nh)$ but
only to $f^{(k)}(0)$, which is given by a weighted sum whose weights
are subsumed under~$c_i$. Similarly, $\beta$ is related to weights
subsumed under~$d_i$.\, So, to our initial concession that ``the
coefficients $c_i$ and $d_i${\tiny\,}\ldots\ may depend on $\alpha$
and~$\beta$,'' we could add ``\emph{respectively}.'' Moreover, the
change-of-variable ${f(t)\!=\!g(u)}$ where ${u\!=\!nh\!-\!t}$ (whence
${t\!=\!nh\!-\!u}$ and ${dt\!=\!-du}$) transforms rule~\eqref{e-cd}
into
\begin{equation}\label{e-cdswap}
  \int_{-\beta h}^{(n\text{+}\alpha)h} \!\!g(u)\,du \,\approx~
    h\sum_{i\text{=0}}^n g(ih) \,+\,
    h\sum_{i\text{=0}}^m d_i \text{~\!}g(ih) \,+\,
    h\sum_{i\text{=0}}^m c_i \text{~\!}g\big((n\!-\!i)h\big) \,,
\end{equation}
which is the same rule except that $\alpha$ and~$c_i$ have swapped
places with $\beta$ and~$d_i$. So if the rule is consistent, the
corrections $d_i$ depend on~$\beta$ as the corrections $c_i$ depend
on~$\alpha$. In~particular, if~${\beta\!=\!\alpha}$, then
${d_i\!=\!c_i}${\tiny\,} and rule~\eqref{e-cd} reduces~to
\begin{equation}\label{e-cdsym}
  \int_{-\alpha h}^{(n\text{+}\alpha)h} \!\!f(t)\,dt \,\approx~
    h\sum_{i\text{=0}}^n f(ih) \,+\,
    h\sum_{i\text{=0}}^m c_i\Big[f(ih) + f\big((n\!-\!i)h\big)\Big] \,.
\end{equation}
This special case, by its symmetry about ${t\!=\!nh\text{/2}}$,\, is
exact if $f(t)$ is any \emph{odd} power of ${(t-nh\text{/2})}${\tiny{}
}\textemdash so~that, if it is exact for a polynomial of \emph{even}
degree~$m$, it is also exact for degree~$m\text{+1}$. This raising of
the maximum degree for exactness does not happen when
${\beta\!\ne\!\alpha}$. And when it \emph{does} happen
(when~${\beta\!=\!\alpha}$), it does not change the order of the error
for a general analytic\,$f$;\, it~happens because when $f$ is of
degree~$m\text{+1}$, the coefficients of $h^{m\text+2}$
in~\eqref{e-cdsym} are matched by the antisymmetry of the
highest-power term in\,$f$, whereas a more general\,$f$ generally
breaks the antisymmetry.\, More precise error bounds for the case
${\alpha\!=\!\beta\!=\!0}$ are given by Barrett~\cite{barrett-52},
Martensen~\cite{martensen-70}, and De~Swardt \& De\,Villiers
\cite[p.~\!131]{deSwardt-deVilliers-00}
(citing~\cite[pp.\,161--3]{martensen-64}).

\section{Finding the solution}

Given that there exist coefficients $c_i$ depending on $\alpha$ alone,
and $d_i$ depending identically on $\beta$ alone, which make
rule~\eqref{e-cd} correct in a certain sense for all functions $f(t)$
of a certain class, we~can choose any convenient member of that class
for the purpose of finding the coefficients. But what is
``convenient''?

First hint: If the rule works for arbitrary~$n$, it must work as
${n\!\to\!\infty}$, provided of course that the integral converges, in
which case the integrand and hence the right-hand sum in~\eqref{e-cd}
must go to zero, so~that we are left with
\begin{equation}\label{e-c}
  \int_{-\alpha h}^{\,\infty} \!f(t)\,dt \,\approx~
  h\sum_{i\text{=0}}^\infty f(ih) \,+\, h\sum_{i\text{=0}}^m c_i ~\!f(ih) \,,
\end{equation}
in which both sides are functions of~$h$. In~the Taylor expansions of
the two sides about ${h\!=\!0}$, the constant terms automatically
match because both sides of~\eqref{e-c} approach the same integral as
${h\!\to\!0}$. And by adjusting the $m\text{+1}$ coefficients~$c_i$,
we~should be able to match the terms in $h^1$ to~$h^{m\text{+1}}$,
so~that the error is~$O(h^{m\text{+2}})$.

The rest of the argument takes copious hints from
Fornberg~\cite{fornberg-11,fornberg-21}, who took less-copious hints
from Fr\"{o}berg~\cite[pp.\,194--6]{froeberg-65} and a fragment of
a~letter by \textbf{James~Gregory} to John Collins
\cite[at~pp.\,208--9]{gregory-1670} dated~\textbf{1670}\textemdash the
year before Newton stated the ``3/8~Simpson rule''\!, and more than
40~years before Cotes computed closed ``Newton--Cotes'' weights for up
to 11~points \cite[p.~\!130]{hairer-wanner-08}. My~generalization via
the parameter~$\alpha$ is largely anticipated by Fornberg \&
Lawrence~\cite[pp.\,4--5]{fornberg-lawrence-23}, whose parameter~$\xi$
corresponds to my~$-\alpha$. Their approach is less general in that
they restrict the range of the parameter (because they are interested
in dealing with discontinuities between samples), but \emph{more}
general in that they use some degrees of freedom to reduce
oscillations in the weights. (And it is more detailed than mine in
some ways, as noted below.)

Second hint: Newton's method of polynomial interpolation
\cite[pp.\,10--12]{hairer-wanner-08} suggests that the
coefficient-matching can be simplified by rewriting~\eqref{e-c} in
terms of \emph{differences} instead of ordinates. If we define the
operator $\Delta$ by\vspace{-1ex}
\begin{equation}\label{e-Dh}
  \Delta f(t) = f(t\text{+}h)-f(t) \vspace{-1ex}
\end{equation}
so that
\begin{equation}\label{e-Dhi}\begin{split}
  \Delta^0 f(0) &= f(0) \\
  \Delta^1 f(0) &= f(h) - f(0) \\
  \Delta^2 f(0) &= f(2h) - 2f(h) + f(0) \\
  \vdots~~ & \\[-.5ex]
  \Delta^m f(0) &= \sum_{j\text{=0}}^m \textstyle\binom{m}{j}(-1)^{m-j} f(jh) ~,
\end{split}\end{equation}
then the second sum in \eqref{e-c}, namely
\begin{equation}\label{e-csum}
  \sum_{i\text{=0}}^m c_i ~\!f(ih) \,, \vspace{-1ex}
\end{equation}
can be written in the form
\begin{equation}\label{e-asum}
  \sum_{k\text{=0}}^m b_k ~\!\Delta^k f(0) \,.
\end{equation}
For, if we equate the last two sums and expand the latter with the
aid of~\eqref{e-Dhi}, we~get\vspace{-1ex}
\begin{equation}
  \sum_{k\text{=0}}^m
    \Bigg(b_k \sum_{j\text{=0}}^k {\textstyle\binom{k}{j}}(-1)^{k-j} f(jh)\Bigg)
  ~=\, \sum_{i\text{=0}}^m c_i ~\!f(ih)
\end{equation}
or, equating the coefficients of $f(ih)$,
\begin{equation}\label{e-ci}
  \sum_{k\text{=}i}^m \textstyle\binom{k}{i}(-1)^{k-i} \text{~\!}b_k \,=~ c_i
  \text{~~}; \quad i\in[0\,\text{..}\,m] \,.
\end{equation}
This is an upper-triangular unit-diagonal system of linear equations,
which can be solved for $b_m$ to~$b_0$ (in~that order) by
back-substitution (Fornberg \&
Lawrence~\cite[p.~\!3]{fornberg-lawrence-23} show the equations in
matrix form). Thus, given the corrections~$c_i$, we~can find the
difference coefficients~$b_k${\tiny{} }\textemdash as claimed. [And of
  course, given the~$b_k$, we~can use~\eqref{e-ci} to find
  the~$c_i$.]\, Substituting \eqref{e-asum} for~\eqref{e-csum}
in~\eqref{e-c}, we~obtain the rule in the desired form
\begin{equation}\label{e-b}
  \int_{-\alpha h}^{\,\infty} \!f(t)\,dt \,\approx~
  h\sum_{i\text{=0}}^\infty f(ih) \,+\, h\sum_{k\text{=0}}^m b_k ~\!\Delta^k f(0) \,,
\end{equation}
and the problem is to find the $m\text{+1}$ constants $b_k$ which
equate the coefficients of the powers of~$h$ from $h^1$
to~$h^{m\text{+1}}$\!.

So a ``convenient'' choice of $f(t)$ should turn the first sum
in~\eqref{e-b} into something tractable\textemdash e.g., a decaying
geometric series. If we choose
\begin{equation}
  f(t) = e^{-st\text{/}h} \vspace{-1ex}
\end{equation}
where $\text{Re}(s)\!>\!0$, then\vspace{-1ex}
\begin{equation}
  \sum_{i\text{=0}}^\infty f(ih)
  \,= \sum_{i\text{=0}}^\infty \big(e^{-s}\big)^i = \frac{1}{1-e^{-s}} \vspace{-1ex}
\end{equation}
and
\begin{equation}
  \int_{-\alpha h}^{\,\infty} \!f(t)\,dt
  \,=\, \frac{e^{-st\text{/}h}}{-s\text{/}h}~\!\Bigg|_{t=-\alpha h}^{t\to\infty}
  =\, h\,\frac{e^{\alpha s}}{s}
\end{equation}
and
\begin{equation}
  f(0) = 1 \,,
\end{equation}
and increasing $t$ by~$h${\tiny\,} multiplies $f(t)$ by~$e^{-s}$
so~that, in operational terms,\vspace{-.5ex}
\begin{equation}
  \Delta = \big(e^{-s} - 1\big) \,.
\end{equation}
If we make these four substitutions in~\eqref{e-b}, we~can cancel~$h$
and obtain
\begin{equation}\label{e-s}
  \frac{e^{\alpha s}}{s} \,\approx\,
  \frac{-1}{e^{-s} \!-\! 1} \,+ \sum_{k\text{=0}}^m b_k\big(e^{-s} \!-\! 1\big)^k .
  \vspace{-1ex}
\end{equation}
Now if we put
\begin{equation}
  x = e^{-s} \!-\! 1 \,,
\end{equation}
so that ${e^{\alpha s}\!=\!(\text{1+}x)^{-\alpha}}$,\,
${s\!=\!-\ln(\text{1+}x)}$,\, and\, ${x\!\to\!0^-}$ as
${s\!\to\!0^+}$, then \eqref{e-s} becomes
\begin{equation}
  \frac{-(\text{1+}x)^{-\alpha}}{\ln(\text{1+}x)} \,\approx\,
  \frac{-1}{x} \,+ \sum_{k\text{=0}}^m b_k ~\!x^k , \vspace{-1ex}
\end{equation}
i.e.,
\begin{equation}\label{e-x}
  \Bigg(\sum_{k\text{=0}}^m b_k ~\!x^k\Bigg) \,\ln(\text{1+}x)
  \,\approx~ \frac{\ln(\text{1+}x)}{x} - (\text{1+}x)^{-\alpha} .
\end{equation}
Taking the geometric series for $(\text{1+}x)^{-1}$ and integrating
term-by-term (putting ${x\!=\!0}${\tiny\,} to set the constant),
we~get\vspace{-1ex}
\begin{equation}
  \ln(\text{1+}x) \,=
  \sum_{i\text{=0}}^\infty \frac{(-1)^i}{i\text{+1}}~\!x^{i\text{+1}} \vspace{-.5ex}
\end{equation}
which, upon dividing by $x${\tiny\,} and renaming the counter,
becomes\vspace{-.5ex}
\begin{equation}
  \frac{\ln(\text{1+}x)}{x}
  \,= \sum_{j\text{=0}}^\infty \frac{(-1)^j}{j\text{+1}}~\!x^j . \vspace{-.5ex}
\end{equation}
The remaining term in~\eqref{e-x} has the binomial expansion\vspace{-1ex}
\begin{equation}
  (\text{1+}x)^{-\alpha}
  = \sum_{j\text{=0}}^\infty \textstyle\binom{-\alpha}{j} ~\!x^j . \vspace{-.5ex}
\end{equation}
With these three substitutions, equation~\eqref{e-x} becomes\vspace{-.5ex}
\begin{equation}\label{e-xser}
  \Bigg(\sum_{k\text{=0}}^m b_k ~\!x^k\Bigg)
  \Bigg(\sum_{i\text{=0}}^\infty \frac{(-1)^i}{i\text{+1}}~\!x^{i\text{+1}}\Bigg)
  ~\approx\,
  \sum_{j\text{=0}}^\infty\bigg[
   \frac{(-1)^j}{j\text{+1}} - {\textstyle\binom{-\alpha}{j}}
  \bigg] \,x^j . \vspace{-.5ex}
\end{equation}

Now we can equate coefficients in~\eqref{e-xser}. On the left side,
there is no term in~$x^0$, due to the index~$i\text{+1}$. On~the
right, the coefficient of~$x^0$ is\vspace{-.5ex}
\begin{equation}
  \frac{(-1)^0}{\text{0+1}} - \textstyle\binom{-\alpha}{0}
  = 1\!-\!1 = 0 \,,
\end{equation}
which agrees with the left side. So the coefficients~$b_k$ must be
fixed so as to match the coefficients of $x^1$
to~$x^{m\text{+1}}$\!.\, If~the product-of-sums on the left is
expanded, there will be $j$~terms in~$x^j$, with $k$ ranging from 0
to~${j\!-\!1}$, and $i$ ranging from ${j\!-\!1}$ to~0
respectively. So, to equate the coefficients of~$x^j$, we take the
second sum on the left inside the first sum, select the inner term
with ${i\!=\!j\!-\!k\!-\!1}$, and take the outer sum up to
${k\!=\!j\!-\!1}$, obtaining\vspace{-.5ex}
\begin{equation}
  \sum_{k\text{=0}}^{j-1} \frac{(-1)^{j-k-1}}{j\!-\!k} ~\!b_k
  \,=~ \frac{(-1)^j}{j\text{+1}} - \textstyle\binom{-\alpha}{j}
  \text{~~}; \quad j\in[1\,\text{..}\,m\text{+1}] \,. \vspace{-.5ex}
\end{equation}
To minimize confusion, let us rename the dummy index\,$j$ as
$i\text{+1}$, so~that both indices count from zero. This
yields\vspace{-.5ex}
\begin{equation}\label{e-bk}
  \sum_{k\text{=0}}^{i} \frac{(-1)^{i-k}}{i\!-\!k\text{+1}} ~\!b_k
  \,=~ \frac{(-1)^{i\text{+1}}}{i\text{+2}} - \binom{-\alpha}{~\!i\text{+1}}
  \text{~~}; \quad i\in[0\,\text{..}\,m] \,,
\end{equation}
in which the coefficient of~$b_k$ is 1 if~${k\!=\!i}$, and there are
no terms for~${k\!>\!i}$.

So~\eqref{e-bk} is a \emph{lower}-triangular unit-diagonal system of
linear equations in~$b_k$, which can be solved for $b_0$ to~$b_m$
(in~that order) by \emph{forward} substitution. This forward order
means that we can increase~$m$, adding more equations, without
invalidating the solutions found so~far. But, having found as many
coefficients~$b_k$ as we want, we~then need to find the
corrections~$c_i$ by direct substitution into the
\emph{upper}-triangular system~\eqref{e-ci}, in which higher-index
values of~$b_k${\tiny\,} \emph{do} affect lower-index values of~$c_i$.

Fornberg \& Lawrence \cite[p.~\!4]{fornberg-lawrence-23} give explicit
formulae showing the variation of $b_k$ with their parameter~$\xi$
(our~$-\alpha$) for selected~$k$, and the asymptotic behavior of $b_k$
for extreme values of~$\xi$, the latter behavior being relevant to the
quest for high-order rules with many well-behaved non-unit
weights. The following survey, in~contrast, ignores their restriction
on the parameter (${0\!\le\!\alpha\!<\!1}$ in our notation) and seeks
rules with relatively \emph{few} non-unit weights.

\section{Examples}

As Fornberg suggests~\cite{fornberg-11}, we might reasonably solve the
equations in {\small MATLAB} if we want the coefficients in decimal
form, or in Wolfram Mathematica if we want them in exact rational
form.\, I~used a spreadsheet! Values of~$k$ were filled across the
top, and values of~$i$ down the left-hand edge. For the binomial
coefficient in~\eqref{e-bk}, the value of $-\alpha$ was entered
manually into the appropriate cell\textemdash the most consequential
cell in the sheet\textemdash and subsequent values were built
recursively. The matrix-inversion and matrix-multiplication functions
were used where convenient. In~the course of the inquiry, I~computed
two columns of corrections~$c_i$ for which exact rational values were
desirable but not obvious; for these I~found a common denominator by
expanding a decimal value in a continued fraction.

\subsection{Validation: Reproducing known rules}

\paragraph{Case \boldmath$\alpha\!=\!0${\tiny\,}:} This is the case
considered by Fornberg, after Gregory, assuming \textit{ab~initio}
that the outermost nodes coincide with the terminals. In~his
``Table~1'' \cite[p.~\!170]{fornberg-21}, where his $p$ is
our~$m\text{+2}$, Fornberg gives the corrections, which are duly
reproduced by our eqs.\,\eqref{e-bk}
and~\eqref{e-ci}.\footnote{\,\small My spreadsheet went as far as
${m\!=\!5}$.}\, E.g., the corrections for ${m\!=\!2}$ are\smallskip

$-\tfrac{\,5\,}{\,8\,}\,,\, \tfrac{\,1\,}{\,6\,}\,,\, -\tfrac{1}{24}\,,$\smallskip

\noindent which, when applied from each end of a unit-weight rule with
six or more points, give the ``Lacroix~rule''
\cite[p.~\!131]{deSwardt-deVilliers-00}{\tiny\ }\textemdash a
Gregorian rule with its own name, having the same order of accuracy as
the ``Simpson'' rules. And if we apply the four corrections for
${m\!=\!3}$ from each end of a sufficiently long unit-weight rule, and
then (say) halve~$h$ and double~$n$, we~find that the error is
$O(h^5)$, which is one better than the ``Simpson''
rules~\cite{putland-25-greg}; and so on.

Moreover, Hamming~\cite[pp.\,342--4]{hamming-73} notes that if we
apply a Gregorian correction sequence from each end of the unit-weight
rule of the \emph{same length} (i.e., if ${n\!=\!m}$), we~get the
standard closed Newton--Cotes rule of that length: ${m\!=\!1}$ gives
the trapezoidal rule,\, ${m\!=\!2}$ gives ``Simpson's 1/3 rule''\!,\,
${m\!=\!3}$ gives ``Simpson's 3/8 rule''\!,\, ${m\!=\!4}$ gives
``Boole's rule''\!,\, etc. This, he~says, is ``perhaps the simplest
way to find the actual coefficients'' of the N--C
rules~\cite[pp.\,342]{hamming-73}.

We should add that if $m$ is \emph{even} (so that the number of
corrections is odd), we~get a standard closed N--C rule not only by
applying the corrections to $m\text{+1}$ unit weights, but also by
applying them to $m\text{+2}$ unit weights. Thus the trapezoidal rule
(with two points) is obtained from ${m\!=\!1}$ (two corrections) or
${m\!=\!0}$ (one correction), and ``Simpson's 3/8 rule'' (with four
points) is obtained from ${m\!=\!3}$ (four corrections) or ${m\!=\!2}$
(three corrections), and so~on.

But if $m$ is \emph{odd} (so that the number of corrections is even),
we~get a standard closed N--C rule only by applying the corrections to
$m\text{+1}$ unit weights, \emph{not} by applying them to $m\text{+2}$
unit weights. Thus we~do not get ``Simpson's 1/3 rule'' (three points)
from ${m\!=\!1}$ (two corrections), nor ``Boole's rule'' (five points)
from ${m\!=\!3}$ (four corrections), although the associated Gregorian
rules still have the full expected accuracy, with
error~$O(h^{m\text{+2}})$.

The Gregorian rule for ${m\!=\!1}$ has an alternative explanation.
The ``corrected'' composite trapezoidal rule, which uses derivatives
at the terminals, is two orders more accurate than the uncorrected one
(that~is, it~has the same order as the ``Simpson'' rules). If~the
Gregorian corrections for ${m\!=\!1}$, namely
${-\tfrac{7}{12}~\!,~\!\tfrac{1}{12}}$, are re-expressed as
corrections to the composite trapezoidal rule, they become
${-\tfrac{1}{12}~\!,~\!\tfrac{1}{12}}$; and this sequence is
recognizable as a finite-difference estimate of the end ``correction''
to the trapezoidal rule,\footnote{\,\small The end corrections are
given by the term in~$h^2$ in eq.\,(2) of
Fornberg~\cite{fornberg-21}.} taking the derivative at the distance
$h\text{/2}$ from the terminal instead of \emph{at} the terminal, and
thereby giving an order of accuracy \emph{between} the uncorrected and
``corrected'' trapezoidal rules.

In general, the corrections for \emph{even}~$m$ make Gregory's rule
exact for degrees up to $m\text{+1}$, like the closed
${(m\text{+1})}$-point and ${(m\text{+2})}$-point N--C rules, which
are the unique closed equispaced rules of their lengths that are exact
up to that degree. And the corrections for \emph{odd}~$m$ make
Gregory's rule exact for degrees up to~$m$, like the closed
${(m\text{+1})}$-point N--C rule, which is the unique closed
equispaced rule of its length that is exact up to that degree.  Thus
Gregory's method must generate every closed Newton--Cotes
rule\textemdash \emph{twice} if the rule has an even number of points
(for an odd-degree interpolating polynomial).

\smallskip

Yet Gregory's letter to Collins~\cite{gregory-1670} predates every
closed N--C rule except the trapezoidal rule and Kepler's barrel rule
(also known as Simpson's 1/3 rule). Five years after he wrote that
letter, Gregory was dead. Tradition holds that he suffered a stroke
while showing his students the moons of Jupiter, whereas the earliest
surviving account says: ``By a cold catched in the castle, he grew
blind in on[e] night, and shortly after
dyed''~\cite{frarghall-1676}. He was 36.

{\setlength{\parskip}{1ex minus .2ex}

\vspace{-2ex}

\paragraph{Case \boldmath$\alpha\!=\!1${\tiny\,}:} This case departs
from Gregory/Fornberg by yielding \emph{open} rules whose outermost
nodes are one step in from the terminals. By~analogy with the
preceding case, the corrections given by eqs.\,\eqref{e-bk}
and~\eqref{e-ci} for \emph{even}~$m$ should yield the open
${(m\text{+1})}$-point and ${(m\text{+2})}$-point N--C rules (listed
by Weisstein~\cite{weisstein-25} for up to 7~points), whereas the
corrections for \emph{odd}~$m$ should yield the open
${(m\text{+1})}$-point N--C rule.\, Let us check.

\noindent For~${m\!=\!0}$, the sole correction is
${c_0\!=\!\tfrac{\,1\,}{\,2\,}}$. When applied (twice) to the 1-point
unit-weight rule, this gives the single weight~2, which agrees with
the 1-point open N--C rule. Applied from each end of the 2-point
unit-weight rule, it~gives the weights of the 2-point open N--C rule:

$\tfrac{\,3\,}{\,2\,}\,,\, \tfrac{\,3\,}{\,2\,}\,.$

\noindent For~${m\!=\!1}$, the correction sequence is

$\tfrac{11}{12}\,,\, -\tfrac{5}{12}\,.$

\noindent Applied from each end of the 2-point unit-weight rule, this
gives the 2-point open N--C rule again.

\noindent For~${m\!=\!2}$, the correction sequence is

$\tfrac{31}{24}\,,\, -\tfrac{\,7\,}{\,6\,}\,,\, \tfrac{\,3\,}{\,8\,}\,.$

\noindent Applied from each end of the 3-point unit-weight rule, this
gives the weight sequence

$\tfrac{\,8\,}{\,3\,}\,,\, -\tfrac{\,4\,}{\,3\,}\,,\, \tfrac{\,8\,}{\,3\,}\,,$

\noindent which is the 3-point open N--C rule. And applied from each
end of the 4-point unit-weight rule, it~gives the weight sequence

$\tfrac{55}{24}\,,\, \tfrac{5}{24}\,,\, \tfrac{5}{24}\,,\, \tfrac{55}{24}\,,$

\noindent which is the 4-point open N--C rule.

\noindent For~${m\!=\!3}$, the correction sequence is

$\tfrac{1181}{720}\,,\, -\tfrac{1593}{720}\,,\, \tfrac{1023}{720}\,,\, -\tfrac{251}{720}\,.$

\noindent Applied from each end of the 4-point unit-weight rule, this
gives the 4-point open N--C rule again.

\noindent For~${m\!=\!4}$, the correction sequence is

$\tfrac{2837}{1440}\,,\, -\tfrac{5086}{1440}\,,\, \tfrac{4896}{1440}\,,\, -\tfrac{2402}{1440}\,,\, \tfrac{475}{1440}\,.$

\noindent Applied from each end of the 5-point unit-weight rule, this
gives the simple but oscillatory weight sequence

$\tfrac{33}{10}\,,\, -\tfrac{21}{5}\,,\, \tfrac{39}{5}\,,\, -\tfrac{21}{5}\,,\, \tfrac{33}{10}\,,$

\noindent which is the 5-point open N--C rule.~ So far: so good.

\vspace{-2ex}

\paragraph{Case \boldmath$\alpha\!=\!-1${\tiny\,}:} This yields
rules for which the outermost nodes are one step \emph{outside} the
range of integration.

\noindent One of these rules is easily confirmed. For~${m\!=\!2}$ the
computed corrections are

$-\tfrac{25}{24}\,,\, -\tfrac{\,1\,}{\,2\,}\,,\, \tfrac{1}{24}\,.$

\noindent When these are applied to a sufficiently long sequence of
unit weights, the first three weights are

$-\tfrac{1}{24}\,,\, \tfrac{\,1\,}{\,2\,}\,,\, \tfrac{25}{24}\,.$

\noindent The respective weights for the composite trapezoidal rule
(with the range of integration beginning at the \emph{second} node)
are

$0\,,\, \tfrac{\,1\,}{\,2\,}\,,\, 1$

\noindent so that, by subtraction, the corrections to the composite
\emph{trapezoidal} rule given by~${m\!=\!2}$ are

$-\tfrac{1}{24}\,,\, 0\,,\, \tfrac{1}{24}\,.$

\noindent The corresponding contribution to the right-hand side
of~\eqref{e-cd} is
\begin{equation}\label{e-dctr}
  h\Big[{-}\tfrac{1}{24} f(0) + \tfrac{1}{24} f(2h)\Big]
  \,=\, \tfrac{1}{12} ~\!h^2 \,\frac{f(2h)-f(0)}{2h}
  \,\approx\, \tfrac{1}{12} ~\!h^2 f'(h) \,,
\end{equation}
where the right-hand expression is the standard left-hand correction
in the ``corrected trapezoidal rule'' (the argument $h$ is the lower
limit of integration). Thus, by taking ${\alpha\!=\!-1}$ and
${m\!=\!2}$, we~get a discretized corrected composite trapezoidal
rule. An equivalent rule is given by Weisstein~\cite{weisstein-25},
who describes it as a ``2-point open extended formula'' without
further explanation. For a single interval (single step), this rule
has the weights

$-\tfrac{1}{24}\,,\, \tfrac{13}{24}\,,\, \tfrac{13}{24}\,,\, -\tfrac{1}{24}\,,$

\noindent which we have obtained by setting ${\alpha\!=\!-1}$,\,
${m\!=\!2}$, and ${n\!=\!3}$. But they can also be obtained by setting
${\alpha\!=\!0}$,\, ${m\!=\!2}$, and ${n\!=\!1}$, so~that
${m\!>\!n}$;\, in~the latter case, the corrections that overshoot the
unit weights are added to 0 instead of~1.

We could pursue higher-order discrete corrections to the trapezoidal
rule by taking ${\alpha\!=\!-2}$ and ${m\!=\!4}$;\, ${\alpha\!=\!-3}$
and ${m\!=\!6}$;\, etc. But, having come this far in order to
demonstrate the effectiveness of our method, let us now use it to
derive some less familiar rules.

\subsection{Application: Correcting the midpoint rule}

\paragraph{Case \boldmath$\alpha\!=\!\text{\bf 1/2}${\tiny\,}:} This
yields open rules for which the outermost nodes are a \emph{half}-step
in from the terminals\textemdash as in the composite midpoint rule,
which of~course is a unit-weight rule, so that the corrections to the
unit-weight rule can also be called corrections to the midpoint rule.
For~${m\!=\!0}$, the sole ``correction'' is ${c_0\!=\!0}$, leaving the
midpoint rule uncorrected. For~${m\!=\!1}$, the order improves by~1.\,
For~${m\!\in\!\big\{2,3,4\big\}}$, we~get rules that we might actually
want to~use. As~\emph{open} rules, they avoid evaluating the integrand
at the terminals, where it may not be defined. Even if the integrand
has a finite limit as we approach the terminal, there is some
convenience in not having to deal with a singularity, or the
possibility of a singularity, \emph{at}~the terminal, wherefore one
might say that open rules are better than closed rules as
\emph{general-purpose} rules.

\noindent For $\alpha\!=\!\text{1/2}$, the corrections for the
nominated values of $m$ are\footnote{\,\small These were posted
in~\cite{putland-25-mid} for ${m\!=\!2}$, and in a comment thereto for
${m\!=\!3}$ and ${m\!=\!4}$.}

$m\!=\!2~\!:~ \tfrac{1}{12}\,,\, -\tfrac{\,1\,}{\,8\,}\,,\, \tfrac{1}{24}\,;$

$m\!=\!3~\!:~ \tfrac{703}{5760}\,,\, -\tfrac{1389}{5760}\,,\, \tfrac{909}{5760}\,,\, -\tfrac{223}{5760}\,;$

$m\!=\!4~\!:~ \tfrac{909}{5760}\,,\, -\tfrac{2213}{5760}\,,\, \tfrac{2145}{5760}\,,\, -\tfrac{1047}{5760}\,,\, \tfrac{206}{5760}\,.$

\noindent The corresponding weights (if there are unit weights left
over) are

$m\!=\!2~\!:~ \tfrac{13}{12}\,,\, \tfrac{\,7\,}{\,8\,}\,,\, \tfrac{25}{24}\,,\, 1\,,\, \ldots\,,\, 1\,,\, \tfrac{25}{24}\,,\, \tfrac{\,7\,}{\,8\,}\,,\, \tfrac{13}{12}\,;$

$m\!=\!3~\!:~ \tfrac{6463}{5760}\,,\, \tfrac{4371}{5760}\,,\, \tfrac{6669}{5760}\,,\, \tfrac{5537}{5760}\,,\, 1\,,\, \text{etc.};$

$m\!=\!4~\!:~ \tfrac{6669}{5760}\,,\, \tfrac{3547}{5760}\,,\, \tfrac{7905}{5760}\,,\, \tfrac{4713}{5760}\,,\, \tfrac{5966}{5760}\,,\, 1\,,\, \text{etc.}$

\noindent The denominator for ${m\!=\!3}$ and ${m\!=\!4}$ was found by
expanding one correction in a continued fraction. (The same approach
to ${m\!=\!5}$ made it clear that the exact rational coefficients
would be unwieldy.)

\noindent As a check, it is worth noting that if we apply the
corrections for ${m\!=\!3}$ from each end of the 4-point unit-weight
rule, we get the same weights\textemdash namely

$\tfrac{13}{12}\,,\, \tfrac{11}{12}\,,\, \tfrac{11}{12}\,,\, \tfrac{13}{12}$

\noindent \textemdash as if we do likewise with the corrections for
${m\!=\!2}$. As~a further check, it~is easily confirmed by experiment
that the resulting rules for ${m\!=\!2}$ and ${m\!=\!3}$ are exact
(to~machine precision) for integrands of degree up to~3 (like the
``Simpson'' rules) while the resulting rule for ${m\!=\!4}$ is exact
for integrands of degree up to~5 (like ``Boole's rule'').

\noindent And as a test, if we integrate ${f(t)\!=\!7t^6}$ from 0
to~1, with 10~nodes, for ${m\!=\!2}$, ${m\!=\!3}$, and ${m\!=\!4}$,
and then double the number of nodes, the errors are reduced by the
approximate factors 13.5, 27.7, and 47.3 respectively, whence it is
not hard to believe that the errors are $O(h^4)$, $O(h^5)$, and
$O(h^6)$ respectively.

\vspace{-2ex}

\paragraph{Case \boldmath$\alpha\!=\!-\text{\bf 1/2}${\tiny\,}:} This
yields rules for which the outermost nodes are a half-step
\emph{outside} the range of integration.

\noindent Again one rule from the series is easily confirmed.
For~${m\!=\!1}$, the computed corrections are

$-\tfrac{23}{24}\,,\, -\tfrac{1}{24}\,.$

\noindent When these are applied to a sufficiently long sequence of
unit weights, the first two weights are

$\tfrac{1}{24}\,,\, \tfrac{23}{24}\,.$

\noindent The respective weights for the composite midpoint rule (the
first midpoint being the \emph{second} node) are

$0\,,\, 1$

\noindent so that, by subtraction, the corrections to that composite
midpoint rule given by~${m\!=\!1}$ are

$\tfrac{1}{24}\,,\, -\tfrac{1}{24}\,.$

\noindent The corresponding contribution to the right-hand side
of~\eqref{e-cd} is
\begin{equation}
  h\Big[\tfrac{1}{24} f(0) - \tfrac{1}{24} f(h)\Big]
  \,=\, -\tfrac{1}{24} ~\!h^2 \,\frac{f(h)-f(0)}{h}
  \,\approx\, -\tfrac{1}{24} ~\!h^2 f'\big(h\text{/2}\big) \,,
\end{equation}
where the right-hand expression is \emph{minus one half} of the
standard left-hand correction in the ``corrected trapezoidal rule''
(the argument $h\text{/2}$ is the lower limit of integration). But it
is clear that $-\text{1/2}$ of the leading-order correction to the
trapezoidal rule is the leading-order correction to the midpoint rule;
e.g., the ``Simpson'' weights
${\big({\tfrac{\,1\,}{\,3\,},\tfrac{\,4\,}{\,3\,},\tfrac{\,1\,}{\,3\,}}\big)}$
are 2/3 of the way from the trapezoidal weights ${(1,0,1)}$ to the
midpoint weights ${(0,2,0)}$. So the rule for
${\alpha\!=\!-\text{1/2}}$ and~${m\!=\!1}$ is a discretized corrected
composite midpoint rule. We~could pursue higher-order discrete
corrections to the midpoint rule by taking ${\alpha\!=\!-\text{3/2}}$
and~${m\!=\!3}$;\, ${\alpha\!=\!-\text{5/2}}$ and~${m\!=\!5}$;\,
etc. (The resulting rules are unusual in that for \emph{odd}~$m$, they
are exact for integrands of degree up to $m\text{+2}$ and thereafter
give errors of~$O(h^{m\text{+3}})$. For they have the same symmetry
about the terminal as the above ``discretized corrected composite
trapezoidal rule''\!, in which, for \emph{even}~$m$, the correction at
the terminal is zero, so that the effective number of corrections is
one \emph{fewer} than would normally be required for the same order of
accuracy.)

\noindent By the way, I~initially derived the rule for
${\alpha\!=\!+\text{1/2}}$ and~${m\!=\!2}$ by treating it as a
corrected midpoint rule, with a different discrete estimate of\,$f'$
at the lower terminal~\cite{putland-25-mid}. But I~used only the
generalized Gregory/Fornberg approach to find the corresponding rules
for ${m\!=\!3}$ and~${m\!=\!4}$.

\subsection{Asymmetrical rules\, \boldmath$(\;\!\beta\!\ne\!\alpha)$}

The examples given so far have used the same value of $\alpha$ at each
end of the range of integration; that~is, in~the notation of
eqs.\,\eqref{e-cd} to~\eqref{e-cdswap}, they have set
${\beta\!=\!\alpha}$. But we can also set $\alpha$ and~$\beta$
independently. This is useful if we have a function sampled at fixed
equispaced abscissae and want to be able to integrate it between
\emph{arbitrary} limits.

\noindent For an illustration, let us take ${\alpha\!=\!\text{1/2}}$
and ${\beta\!=\!0}$, so that the rule is midpoint-like from the left
and closed--N-C--like from the right (such a rule, being open at one
end and closed at the other, is described as \textit{semi-open}).
If~any unit weights remain, the weights for ${m\!=\!2}$ are

$\tfrac{13}{12}\,,\, \tfrac{\,7\,}{\,8\,}\,,\, \tfrac{25}{24}\,,\, 1\,,\, \ldots\,,\, 1\,,\, \tfrac{23}{24}\,,\, \tfrac{\,7\,}{\,6\,}\,,\, \tfrac{\,3\,}{\,8\,}\,,$

\noindent and the weights for ${m\!=\!3}$ are

$\tfrac{6463}{5760}\,,\, \tfrac{4371}{5760}\,,\, \tfrac{6669}{5760}\,,\, \tfrac{5537}{5760}\,,\, 1\,,\, \ldots\,,\, 1\,,\, \tfrac{739}{720}\,,\, \tfrac{211}{240}\,,\, \tfrac{299}{240}\,,\, \tfrac{251}{720}$

\noindent where the last four are given by
Fornberg~\cite[p.~\!8]{fornberg-11}. If we integrate from 0 to~1 with
${h\!=\!\text{2/19}}$ (giving 10~nodes), we~find that the rule for
${m\!=\!2}$ is exact for integrands of degree up to~2 (\emph{not} 3 as
for the symmetrical rules) while the rule for ${m\!=\!3}$ is exact for
integrands of degree up to~3 (as~for the symmetrical rules). For
${f(t)\!=\!5t^4}$, if we reduce $h$ from 2/19 (10~nodes) to~2/39
(20~nodes), the error is reduced by a factor 18.5 for ${m\!=\!2}$,
and~36.4 for ${m\!=\!3}$, whence it is not hard to believe that the
errors are $O(h^4)$ and $O(h^5)$ respectively; recall that, after
eq.\,\eqref{e-cdsym}, the error for a general\,$f$ was expected to be
$O(h^{m\text{+2}})$, regardless of whether the error canceled for $f$
of degree~${m\text{+1}}$ for even~$m$.

\noindent For another illustration, one of the ``single interval
extrapolative rules'' listed by Weisstein~\cite{weisstein-25}, namely
\begin{equation}
  \int_{-h}^0 \!f(t)\,dt \,\approx\,
  h\Big[\tfrac{23}{12}f(0)-\tfrac{\,4\,}{\,3\,}f(h)+\tfrac{5}{12}f(2h)\Big] \,,
\end{equation}
is recognizable as a backward Adams--Bashforth rule, and can be
obtained by setting ${\alpha\!=\!1}$,\, ${\beta\!=\!-2}$, and
${m\!=\!n\!=\!2}$. In~general, \emph{forward} Adams--Bashforth weights
are given by ${\alpha\!=\!-m\!=\!-n}${\tiny\,} and{\tiny\,}
${\beta\!=\!1}$, and forward Adams--Moulton weights by\,
${\alpha=1\!-\!m=1\!-\!n}$\, and{\tiny\,} ${\beta\!=\!0}$.

} 

\section{Conclusion}

The working equations \eqref{e-ci} and~\eqref{e-bk} bear repeating,
and the former bears switching left-to-right for actual use. So,
in~summary, the na\"{i}ve unit-weight equispaced quadrature rule may be
corrected exactly for integrands up to degree~$m$ by adding $m\text{+1}$
``corrections'' to the weights at \emph{each} end, starting with the
outermost weight and working inward. The corrections are given by
\begin{equation}
  c_i \,= \sum_{k\text{=}i}^m \textstyle\binom{k}{i}(-1)^{k-i} \text{~\!}b_k
  \text{~~}; \quad i\in[0\,\text{..}\,m] \,,
\end{equation}
where the $m\text{+1}$ coefficients $b_k$ are the solutions of the
lower-triangular system
\begin{equation}
  \sum_{k\text{=0}}^{i} \frac{(-1)^{i-k}}{i\!-\!k\text{+1}} ~\!b_k
  \,=~ \frac{(-1)^{i\text{+1}}}{i\text{+2}} - \binom{-\alpha}{~\!i\text{+1}}
  \text{~~}; \quad i\in[0\,\text{..}\,m] \,,
\end{equation}
where $\alpha$ is the distance from the limit-of-integration to the
first node, measured inward in step-lengths. The corrections may
overlap, in which case they are cumulative (and if any corrections
overshoot the unit weights, they are added to 0 instead of~1).

The values of $\alpha$ at the two ends need not be the same. If they
\emph{are} the same (``${\beta\!=\!\alpha}$'') and $m$ is \emph{even},
the rule is exact for integrands of degree up to $m\text{+1}$ (instead
of~$m$). Be that as it may, the error in the integral is generally
$O(h^{m\text{+2}})$, where $h$ is the step-length.

Whereas the original purpose of this study was to find corrected
weights for the composite midpoint rule (which, for better or worse,
determined the sign convention for the parameter~$\alpha$), a~wide
variety of closed, open, and extrapolative equispaced rules may be
derived from the same two equations by suitably choosing~$\alpha$,
$\beta$, and $m$, and (for overlapping corrections) the initial number
of unit weights.

\section{Acknowledgment}

I thank Professor Bengt Fornberg for saving me much embarrassment by
giving feedback on an advanced draft of this paper\textemdash and for
his just-published tribute to James Gregory~\cite{fornberg-25}, which
led me directly or indirectly to references
\cite{barrett-52,deSwardt-deVilliers-00,frarghall-1676,hamming-73,martensen-70}.
Any remaining errors or omissions are my own.

\linespread{0.98}\small\raggedright

\markright{References}

\end{document}